\documentclass[10pt]{article}
\usepackage{epsfig,amsmath2000,amsfonts}

\let\mgp=\marginpar \def\marginpar#1{\mgp{\raggedright\tiny #1}}

\def \jmp#1{\marginpar{\%\% \ #1 ---j}}
\let\lbl=\label
\def\label#1{\lbl{#1}\ifinner\else\marginpar{\ref{#1} #1}\ignorespaces\fi}

\def \marginpar#1{}

\def\Bbb{\mathbb}

\def\Ee{{\Bbb E}}   
\def\Et{\Ee^3}    
\def\Z{{\Bbb Z}}    

\def\Msad{{\mathsf M}} 
\def\Scrw{{\mathsf S}} 
\def\Zcrw{{\mathsf Z}} 
\def\Flat{{\mathsf X}} 
\def\Emty{{\mathsf O}} 

\def\pli{\pi}    
\def\plh{\Pi}    
\def\pushed/{pushed-$P_\tau$}
\def\pushedz/{pushed-$P_0$}

\def\vp{\vspace{.5\baselineskip}} \def\vvp{\vspace{.9\baselineskip}}

\def\ddefine#1{{\bf #1}}    
\def\ournum#1{\thesection.\arabic{#1}} 

\newtheorem{thm}{Theorem}[section]
\newtheorem{prop}[thm]{Proposition}
\newtheorem{cor}[thm]{Corollary}
\newtheorem{lem}[thm]{Lemma} 
 
\newtheorem{mainthm}[thm]{Main Theorem}
\def\MainTheorem/{Main Theorem~\ref{mainthm}}

\def\qed{\hfill {\vbox{\hrule\hbox{\vrule height6pt\hskip6pt\vrule}\hrule}}}
\newenvironment{proof}{\smallskip\par\noindent{\bf Proof}\ }{\qed\medskip}

\def\inline#1{${\vcenter{\hbox{\includegraphics{#1}}}}$}
\def\inlinex#1{${\vcenter{\hbox{\includegraphics[scale=.6,angle=90]{#1}}}}$}
\def\inlinexb#1{${\vcenter{\hbox{\includegraphics[scale=.9,angle=90]{#1}}}}$}
\def\inlines#1{${\vcenter{\hbox{\includegraphics[scale=.6]{#1}}}}$}  
\def\inliness#1{$\smash{\vcenter{\hbox{\includegraphics[scale=.6]{#1}}}}$}  
\def\ccaption#1{\caption{#1}} 

\makeatletter
\long\def\@makecaption#1#2{{\narrower\vskip\abovecaptionskip
\sbox\@tempboxa{{\sf #1: #2}}\ifdim\wd\@tempboxa>\hsize {\sf #1: #2\par}
\else\hbox to\hsize{\hfil\box\@tempboxa\hfil}\fi
\vskip\belowcaptionskip}}

\begin{document}\makeatletter

\centerline{\LARGE \sf Cubic Polyhedra}

\vp
\begin{center}\sf\Large
Chaim Goodman-Strauss\footnote{Dept.~Mathematics,
Univ.~Arkansas, Fayetteville AR 72701, {\sl cgstraus@comp.uark.edu}.
Research partially supported by NSF grant DMS-00-72573 and by the
Consejo Nacional de Ciencia y Tecnologia (CONACYT) of Mexico.} \\
John~M.~Sullivan\footnote{Dept.~Mathematics,
Univ.~Illinois, Urbana IL 61801, {\sl jms@math.uiuc.edu}.
Research partially supported by NSF grant DMS-00-71520.}
\end{center}

\vp
\begin{center}\small
Dedicated to W.~Kuperberg on the occasion of his sixtieth birthday,\\
and to the memory of Charles E.~Peck.
\end{center}

\vp
\begin{abstract}
A \emph{cubic polyhedron} is a polyhedral surface whose edges 
are exactly all the edges of the cubic lattice.
Every such polyhedron is a discrete minimal surface,
and it appears that many (but not all) of them
can be relaxed to smooth minimal surfaces
(under an appropriate smoothing flow, keeping their symmetries).
Here we give a complete classification of the cubic polyhedra.
Among these are five new infinite uniform polyhedra and an uncountable
collection of new infinite semi-regular polyhedra.
We also consider the somewhat larger class of
all discrete minimal surfaces in the cubic lattice.
\end{abstract}

\vvp
\section{Introduction}
We define a \ddefine{cubic polyhedron} $P$ to be any polyhedron
whose vertices and edges are exactly the vertices and edges of the cubic
lattice in $\Et$, and which forms an embedded topological surface
(or $2$-manifold) in $\Et$.
(We pick fixed orthonormal coordinates on $\Et$; and we view the cubic
lattice as a cell-complex with vertices at $\Z^{3}$.)
It follows that cubic polyhedra are connected, non-compact and 
unbounded, and of course have faces that are among the square faces of 
the cubic lattice.  Note that we allow adjacent coplanar squares,
and still consider them to be distinct faces of the polyhedron.

A Hamiltonian path in a graph (1-complex) is a connected subcomplex
which forms a 1-manifold and includes every vertex (0-cell) of the graph.
By analogy, we could describe cubic polyhedra as the ``Hamiltonian surfaces''
in the 2-skeleton of the cubic lattice: they are the subcomplexes which
form connected 2-manifolds including every edge of the lattice.
(Banchoff~\cite{Banchoff} used such Hamiltonian surfaces in the 2-skeleton
of the $n$-cube as examples of tight polyhedral surfaces; for $n=6$
his example is a quotient of our cubic polyhedron $P_0$.)

Note that we do not allow the configuration consisting of two touching
cube corners (which would correspond to a pair of opposite triangular
loops in the octahedron) because such a vertex would be a nonmanifold point.

\begin{figure}[ht]
\centerline{\epsfig{file=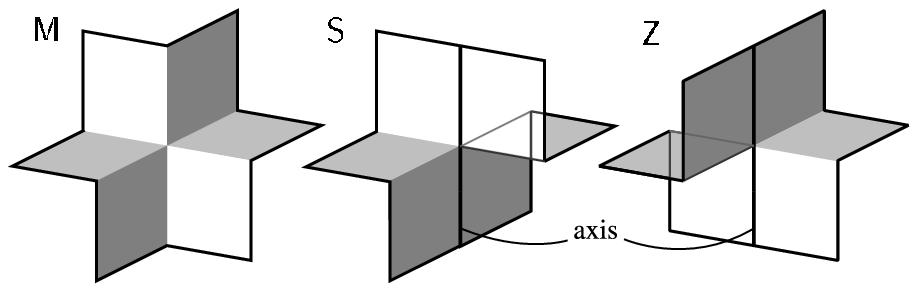}}
\ccaption{Up to rigid motion, there are three possible vertex configurations
in a cubic polyhedron, a monkey-saddle $\Msad$, a left-handed screw $\Scrw$,
and a right-handed screw $\Zcrw$.  The two flanges of a screw vertex lie
along its axis line.}\label{typeMSZ}
\end{figure}

Given a cubic polyhedron $P$, 
every edge in the cubic lattice is incident to two faces of $P$,
and we call it a \ddefine{crease} or a \ddefine{flange} depending on whether
these faces are perpendicular or coplanar.  Similarly, every vertex
has valence six, and we find that (up to isometry) there are only
two possibilities for the configuration of incident squares, as shown
in figure~\ref{typeMSZ}.
To see this, consider the vertex figure of the cubic lattice,
a regular octahedron.  The two vertex configurations correspond
to the two possible Hamiltonian cycles in the 1-skeleton of this octahedron.

In the first possible configuration, a \ddefine{monkey-saddle} ($\Msad$) vertex,
all six incident edges are creases, alternating up and down, making
the vertex indeed a (polyhedral) monkey-saddle.  The normal vector at
an $\Msad$ vertex (meaning the average of the normals to the six
incident faces) points along one of the four body-diagonals;
the monkey-saddle configuration can occur in four possible orientations.

The second configuration, a \ddefine{screw} vertex, has two incident
flanges, in opposite directions along a common \ddefine{axis}.
This configuration comes in left- and right-handed versions, which
we call $\Scrw$ and $\Zcrw$ vertices, respectively.  The normal vector
at a screw vertex points along one of the two face-diagonal lines
perpendicular to the axis.  There are six orientations for an $\Scrw$ vertex
or for a $\Zcrw$ vertex, corresponding to the choice of axis and
normal lines.

Our problem of classifying cubic polyhedra comes down to
figuring out all the ways that these vertex configurations
can be fit together to fill up all of space.   An uncountable number
of polyhedra can be built, with surprising variety; but there
is also certain rigidity which aids our classification,
culminating in our \MainTheorem/.

\vp
\subsection{Basic Constructions}\label{basic}
Before giving some examples, we prove two useful lemmas
which let us extend a finite configuration to a complete polyhedron.

\begin{lem}\label{lem:refl}
Suppose we are given compatible configurations at the vertices
within a rectangular box (which can be finite, infinite or bi-infinite
in each of the coordinate directions).  This can be extended to
a complete cubic polyhedron, by repeated reflection in the sides
of the box (which are planes at half-integer coordinate values).
\end{lem}
\begin{proof}
The edges meeting the boundary of the box do so perpendicularly
(as do their incident faces) no matter whether they are creases or flanges.
Thus the half of each such edge or face within the box reflects
to the half outside; the fact that each vertex within the box
has a legal configuration means the same is true at all reflected
vertices.
\end{proof}

\begin{lem}\label{lem:uniq}
Given the configuration of faces along an edge $e$,
the configuration at either endpoint $v$ of that edge is determined
uniquely by its type $\Msad$, $\Scrw$ or $\Zcrw$.  In particular,
if the vertex configurations at the ends of $e$ are $\Msad$ and $\Msad$,
or $\Scrw$ and $\Zcrw$, then these configurations are mirror images
of each other.
\end{lem}
\begin{proof}
If $e$ is a flange, then its vertices must be screws with axis along $e$,
which (together with the flange normal) fixes their orientation.
If $e$ is a crease, then there are four possible orientations of that crease.
If $v$ is type $\Msad$, it has four possible orientations; if it is type
$\Scrw$ (or $\Zcrw$) it has four possible orientations with axis perpendicular
to $e$.  In any case, the possibilities correspond bijectively to the four
for the crease along $e$.
\end{proof}

These simple lemmas immediately give us the two symmetric
cubic polyhedra shown in figure~\ref{symmpoly}:
\begin{figure}[ht]\centering
 \begin{minipage}[b]{.45\textwidth}\centering\epsfig{figure=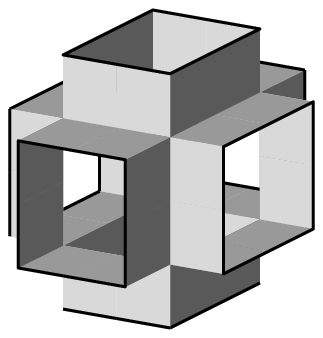}\end{minipage}
 \begin{minipage}[b]{.45\linewidth}\centering\epsfig{figure=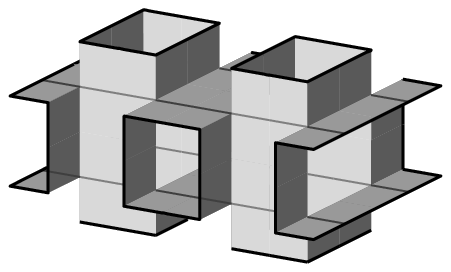}\end{minipage}
 \ccaption{The cubic polyhedra $P_0$ (left) and $P_1$ (right)
are generated by repeated reflection of
a single monkey-saddle or screw vertex, respectively. The thin
black lines in the picture of $P_1$ are its parallel axes.}
\label{symmpoly}
\end{figure}

\begin{prop}
Up to isometry,
there is a unique cubic polyhedron $P_0$ with all monkey-saddle vertices,
and there is a unique cubic polyhedron $P_1$ with all screw vertices whose
types $\Scrw$ and $\Zcrw$ alternate in checkerboard fashion.
\end{prop}
\begin{proof}
Existence follows from Lemma~\ref{lem:refl}, starting with
a single $\Msad$ or $\Scrw$ vertex $v_0$ in a $1\times1\times1$ box.
Rotating the configuration at $v_0$ will just rotate the entire polyhedron.
Uniqueness then follows from Lemma~\ref{lem:uniq}, working outward from $v_0$.
\end{proof}

Of course, $P_0$ is the infinite regular polyhedron $\{4,6\}$
as described by Coxeter~\cite{Cox},
with symmetry group transitive on flags.  We can view $P_1$ as a new
uniform polyhedron:
its faces are all regular polygons and its symmetry group
acts transitively on the vertices, but the faces or edges around
a vertex fall into more than one transitivity class.  Later we will see
four more uniform examples among our cubic polyhedra.

\vp\subsection{Curvature and Topology of Cubic Polyhedra}

Any cubic polyhedron $P$ has, by definition,
six squares meeting at every vertex.
Thus it has an equal quantum $-\pi$ of total Gauss curvature at each vertex.
The surface is an Alexandrov space with curvature bounded above by zero.
This negative curvature spread uniformly throughout space suggests
that $P$ should have nontrivial topology everywhere,
whether or not $P$ is triply periodic.

To examine this, consider an arbitrary loop in the edges of $P$ (which are
the edges of the lattice).  It can be written (homologically) as
a sum of square (four-edge) loops.  Given a square loop $\gamma$,
either $\gamma$ is spanned by a square in $P$ and thus is trivial in $\pi_1(P)$,
or $\gamma$ is a closed geodesic in $P$ and thus (because of the nonpositive
curvature) is nontrivial.  In fact, the square loop $\gamma$ is nontrivial
in $\pi_1(P)$ only if it is nontrivial in $H_1(P)$:
Suppose there is a compact spanning surface $K$ for $\gamma$ within $P$.
If $K$ is just the square convex hull of $\gamma$, then $\gamma$ is nontrivial.
Otherwise $K$ is not contained in the convex hull of $\gamma$,
so it must have an extreme point away from $\gamma$.  But neither
kind of vertex in a cubic polyhedron can be an extreme point.

Incident to any vertex $v$ of the cubic lattice are twelve squares,
four in each of the three directions.  In a cubic polyhedron $P$,
independent of whether $v$ is a monkey-saddle or screw vertex,
exactly two of the four squares in any direction are present in $P$.
The missing squares exhibit nontrivial loops in $H_1(P)$,
which are thus equidistributed in space.  (Note that these loops are
not all independent.)

We note also that every cubic polyhedron $P$ is orientable.
If not, there would necessarily be some orientation-reversing
square loop.  We merely need to check the nine possible $2\times2$
planar diagrams (defined below in Section~\ref{sec:examples})
with a missing central square, to see that this is impossible.

If a cubic polyhedron $P$ has orientation-preserving translational
symmetry with respect to some index $k$ sublattice of the cubic
lattice, then it projects to a compact orientable surface $\overline P$
in the quotient torus (which has volume $k$).  This surface $\overline P$
has $k$ vertices (with total Gauss curvature $-k\pi$),
$3k$ edges and $3k/2$ faces,
so it must have Euler number $-k/2$ and genus $k/4+1$.
(Note, however, that there are $3k/2$ missing squares within this torus.
They form loops that generate $H_1(\overline P)$ but are clearly
not independent since $H_1$ has rank only $k/2+2$.)

For example, $P_0$ and $P_1$ both have translational symmetry with respect to
the even integer lattice $2\Z^3$, with $k=8$, so they have quotients
of genus three.

We can induce a smooth constant-curvature (hyperbolic) metric on
a cubic polyhedron, by giving each
square face the metric of a square in the hyperbolic
plane sized to have internal angle $\pi/3$.

\vp
\subsection{Minimality of cubic polyhedra} 

\begin{figure}[ht]\centering
 \begin{minipage}[b]{.45\linewidth}\centering\epsfig{figure=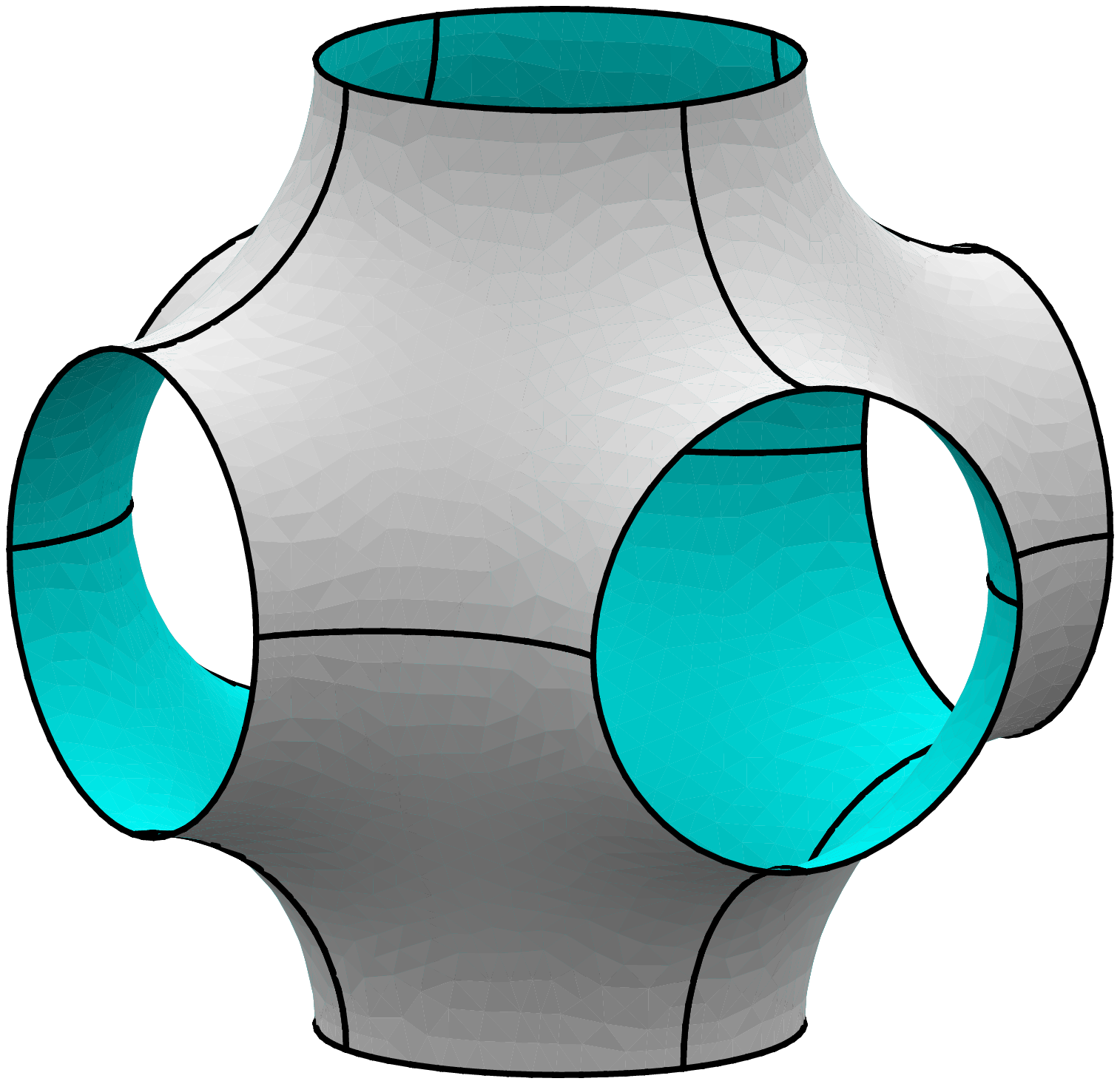,height=1.75in}\end{minipage}
\begin{minipage}[b]{.45\linewidth}\centering\epsfig{figure=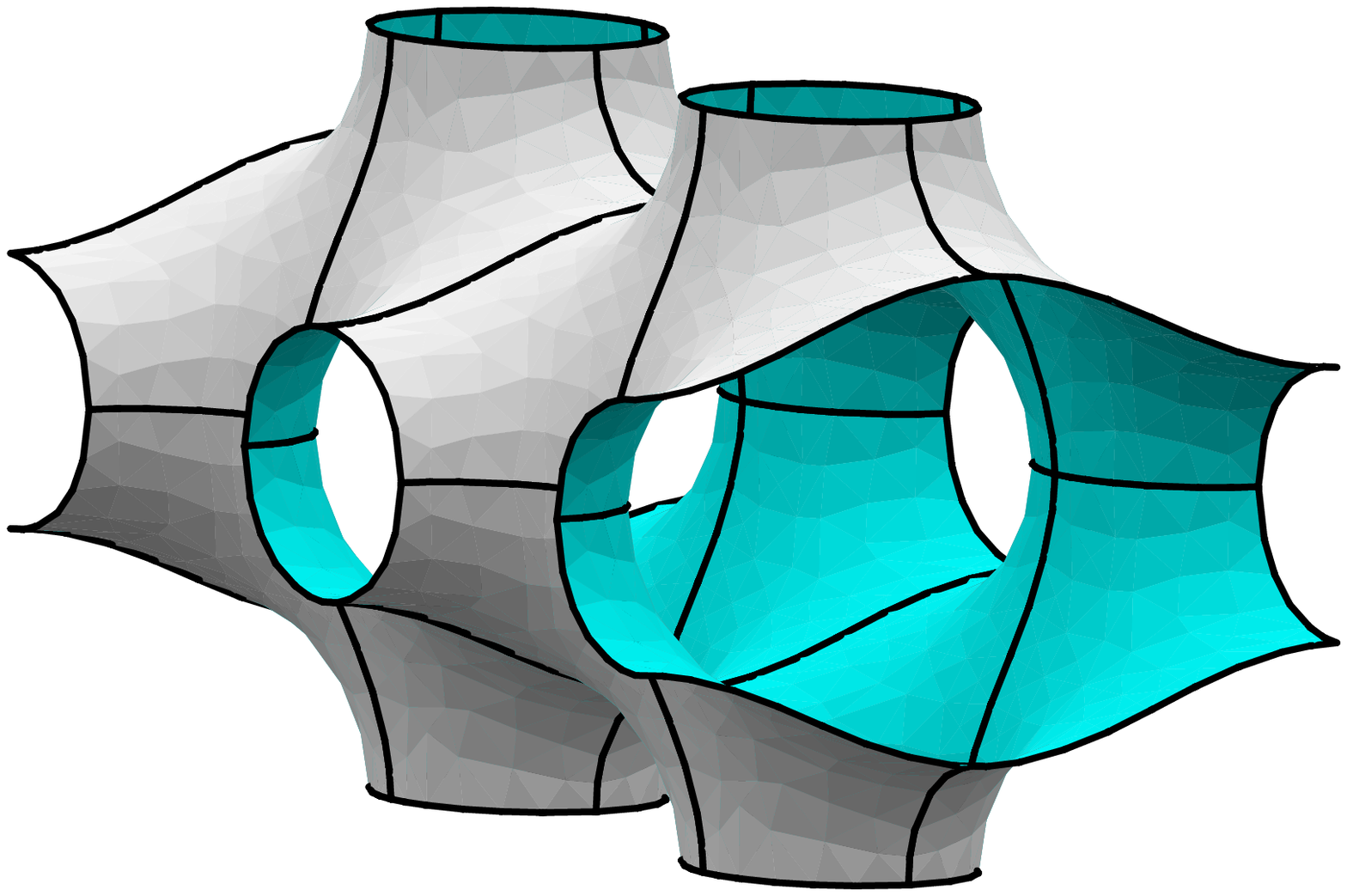,height=1.75in}\end{minipage}
 \ccaption{The polyhedra $P_0$ and $P_1$ relax to triply periodic
minimal surfaces of Schwarz, the P surface (left) and CLP surface (right).
In both cases, the reflection symmetries (across the dark lines) are
preserved under relaxation.}
  \label{sch}
\end{figure}

Our interest in cubic polyhedra arose from the fact that the
first two examples mentioned above were reminiscent of certain
classical triply-periodic minimal surfaces of Schwarz (see~\cite{KP}).
Indeed, we expect that $P_0$ will relax to the P surface
and $P_1$ to the CLP surface.  (See figure~\ref{sch}).

\begin{figure}[ht]
\centering\epsfig{figure=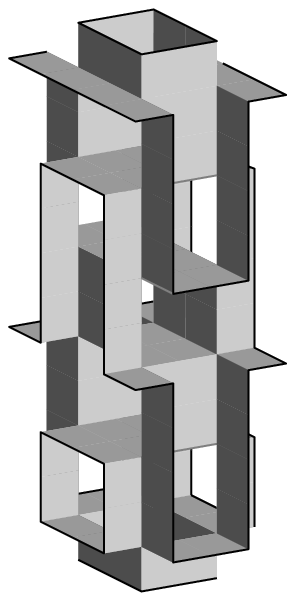}\qquad\epsfig{figure=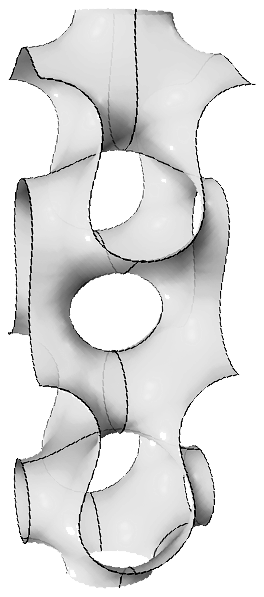,height=2.43in}
\ccaption{A typical example of a cubic polyhedron (left) and
the minimal surface (right) to which it relaxes under the Willmore flow.}
\label{minsurf1}
\end{figure}

We have performed numerical simulations to confirm this, using
Brakke's Evolver~\cite{Brakke}.
We subdivide the faces and let the geometry relax
under a flow which decreases the Willmore bending energy~\cite{HKS}.
(We use this rather than mean-curvature flow, since the minimal
surfaces we flow towards are unstable.)
Further experiments, starting from other cubic polyhedra, indicate
that many, but not all, will relax to minimal surfaces in a similar way. 
(A typical example is shown in figure~\ref{minsurf1}; we have used
pictures like this on posters and
in the first author's \emph{Ptolemy mathcard} series.)
We plan to report on these experiments in a future paper.

Our cubic polyhedra themselves are discrete minimal surfaces,
in the sense of Pinkall and Polthier~\cite{PP}.
Just as a smooth surface is minimal if it
is critical for surface area, a triangulated surface is defined to be
discrete minimal if the first variation of its area is zero, under motion
of any interior vertex.  For a more general polyhedral surface,
we introduce diagonals to triangulate it,
and note that the condition of discrete minimality
is independent of the choice of these diagonals.

Given a surface $P$ made of squares, it is not hard to check that $P$
is discrete minimal at a vertex $v$ exactly when $v$ is the center of mass of
the incident faces or edges.  This happens whenever $v$ has valence six,
as in our cubic polyhedra, but also for exactly one other configuration at $v$:
four coplanar squares.  At the end of this paper, we
briefly consider the calssification of this more general family
of all discrete minimal surfaces within the cubic lattice.

\vp
\subsection{Examples with both types of vertices}\label{sec:examples}

\begin{figure}[htpb]
\centerline{\epsfig{file=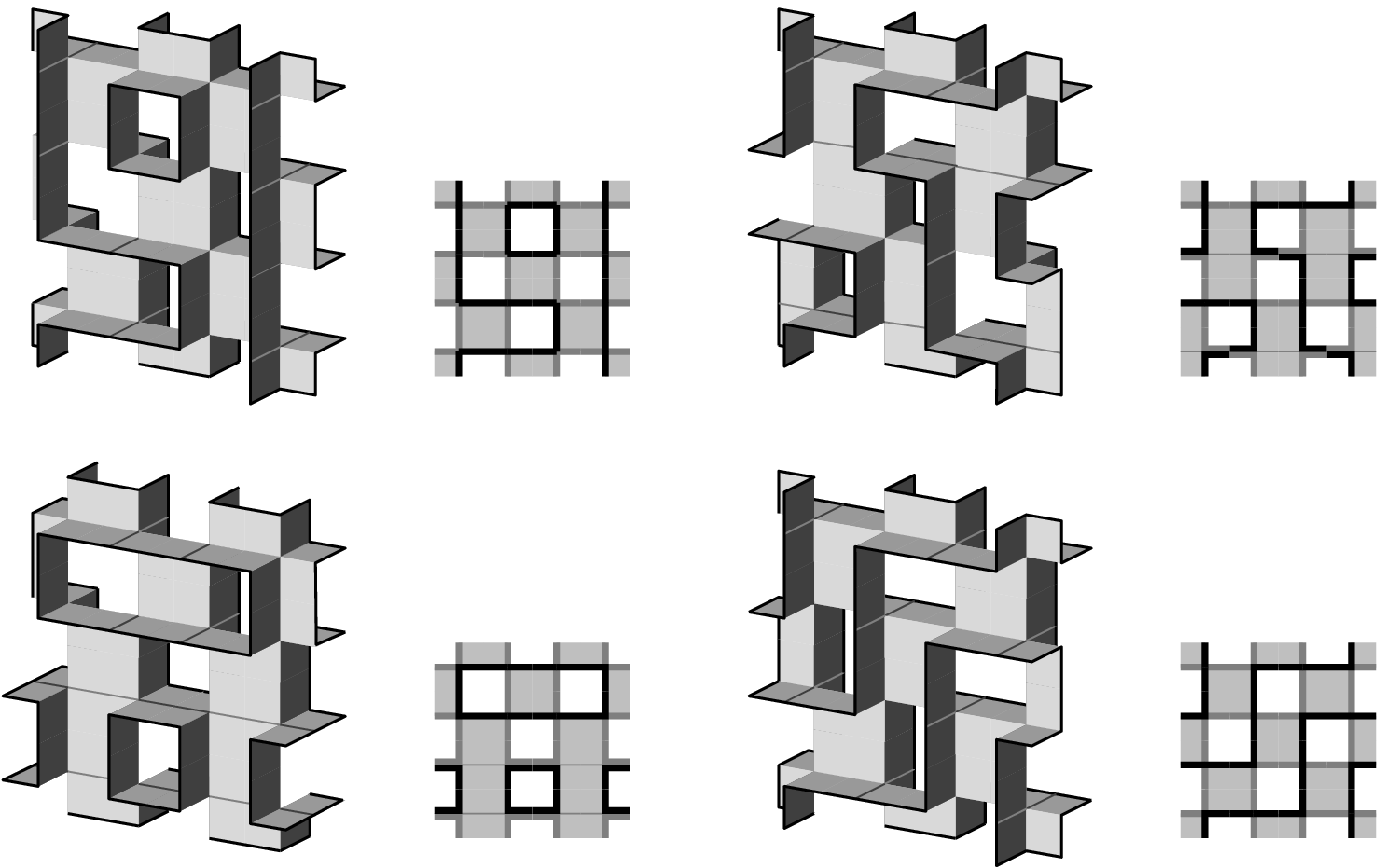}}
\ccaption{Varied examples of cubic polyhedra, with their corresponding planar
diagrams.  These $1\times4\times4$ boxes can be extended by translations in
two directions to $1\times\Z\times\Z$ boxes, and then by reflection
to complete cubic polyhedra.}\label{lotsexamples}
\end{figure}

We can combine monkey-saddle and screw vertices in an astounding variety 
of ways. To generate a cubic polyhedron from any of the complexes in 
figure~\ref{lotsexamples}, we can repeatedly reflect across the front and 
back bounding planes, and translate vertically and to the right and 
left. These examples are meant to suggest that the full class of 
cubic polyhedra is quite large and varied.

We introduce some graphical notation to help discuss examples. 
Given a cubic polyhedron $P$, we use a \ddefine {planar diagram} to 
illustrate a slice through $P$ along some oriented plane $\pli$ in the cubic 
lattice. A planar diagram will have shaded squares corresponding to 
the faces of $P$ in $\pli$, black edges corresponding to the faces of 
$P$ incident to $\pli$ from above, and grey edges corresponding to the faces of 
$P$ incident to $\pli$ from below. So for example,
a monkey-saddle vertex will always appear (up to congruence)
as~\inliness{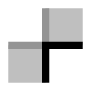}, and a screw vertex will appear (up to congruence)
as either~\inliness{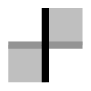} or~\inliness{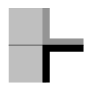},
depending whether or not the axis of the vertex 
is normal to the plane of the diagram.
Figure~\ref{lotsexamples} also shows planar diagrams for the examples.

\vp
\subsection{Configurations around a face in a cubic polyhedron}\label{config}
\begin{figure}[ht]\centering\epsfig{file=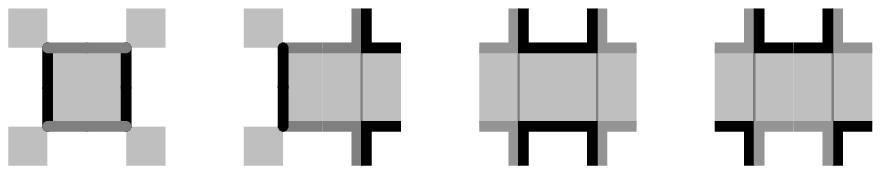}
\ccaption{These are the four possible partial planar diagrams in the
plane of a face $f$. At the left, a normal face has four creases,
alternating up and down.  Next, a face with one flange has three
alternating creases.  Finally, a face with two opposite flanges
also has two opposite creases; there are then exactly two possibilities
for the complete planar diagram around $f$, depending on whether
the creases are to the same side or not.}\label{partialdiag}
\end{figure}

Let $f$ be a (square) face, which we think of as horizontal,
in a cubic polyhedron $P$.  Each of the
four edges of $f$ is either a crease or a flange, but successive
edges cannot both be flanges, because the vertex between them would
then have two perpendicular flanges.  Furthermore, if there are successive
creases, the neighboring squares across them have to alternate up and down
(from the plane of $f$).  We deduce that a (partial) planar diagram in the plane
of $f$ must look like one of the four possibilities in figure~\ref{partialdiag}.

\begin{figure}[ht]\centering
\begin{minipage}[c]{.5\linewidth}\centering\epsfig{file=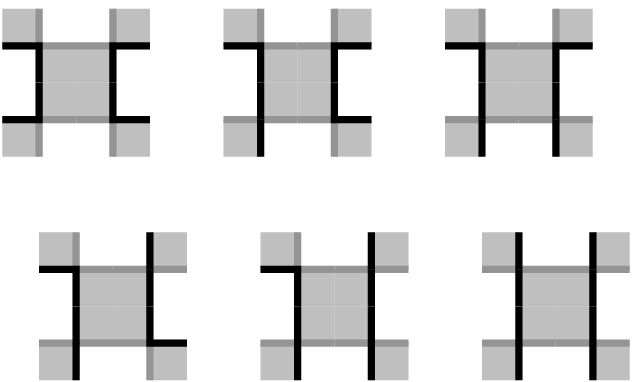}
\ccaption{These six possible planar diagrams around
a normal face, up to isometry and switching the colors of the edges,
are the ways to complete the first partial diagram from
figure~\ref{partialdiag}, corresponding to the six dot diagrams shown
in the text.}\label{towerdiag}\end{minipage}\hfill
\begin{minipage}[c]{.4\linewidth}\centering\epsfig{file=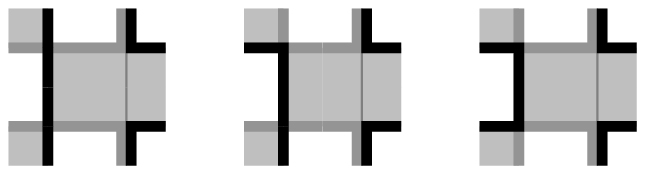, width=2in}
\ccaption{These three possible diagrams around a face with one flange are
the completions of the second partial diagram from
figure~\ref{partialdiag}.}\label{flangediag}\end{minipage}
\end{figure}

We say that $f$ is a \ddefine{normal} face of $P$ if no adjacent faces are
coplanar with $f$, that is, if all its edges are creases.  The four adjacent
faces across these creases must alternate up and down.  (This is the
first case in figure~\ref{partialdiag}.)
The complete configuration of $P$ in a neighborhood of $f$ is then given 
(up to isometry) by one of the six planar diagrams in figure~\ref{towerdiag}.

To see this, consider how monkey-saddle and screw vertices (whose axes
must be vertical) occur in cyclic order around $f$.
There are {\it a priori} six possibilities:
\inline{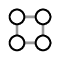}, \inline{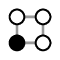}, \inline{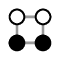},
\inline{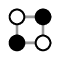}, \inline{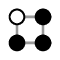}, and \inline{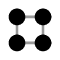},
where these \ddefine{dot diagrams} record the types of vertices
in a slice through $P$ by dots on the integer lattice, with
unfilled dots~\inline{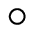} for $\Msad$ vertices and filled dots~\inline{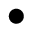}
for screw vertices.  For each of these six possibilities,
the fact that the central square $f$ is normal means that there is
a uniquely determined configuration as in figure~\ref{towerdiag}.

A face $f$ of $P$ which is not normal has, among its edges,
either one flange or two opposite flanges.  If there is
one flange, it is the horizontal axis connecting two screw vertices.
The other two vertices of $f$ can be either monkey-saddles
or screws (with vertical axis); the three resulting configurations
are shown in figure~\ref{flangediag}.

A face $f$ with two (opposite) flanges among its edges must be
in one of the last two configurations shown in figure~\ref{partialdiag}.
Note that the second one, where each flange connects two
screws of the same handedness, is the unique configuration
of a face $f$ with two opposite creases bent to different sides:

\begin{lem}\label{lem:crease}
Suppose a face $f$ in a cubic polyhedron has two opposite edges which
are creases, and the two adjacent faces across these edges
lie on opposite sides of the plane containing $f$.  Then the other
two edges of $f$ are flanges, and each of them is the common axis
of two successive screw vertices of equal handedness.
\qed\end{lem}

We will not illustrate the nine possible diagrams around
a missing square in a cubic polyhedron.  Note, however, that
if all four edges are creases, there are four configurations,
corresponding to the dot diagrams with even numbers of black (or white) dots.
In particular, the diagrams \inline{d2} and \inline{d5}
cannot occur around a missing square.

The local configuration around a normal face $f$ is a $2\times2\times1$
box in one of the six configurations shown in figure~\ref{towerdiag}.
We define a \ddefine{tower} to be one of the (six)
$2\times2\times\Z$ configurations obtained from these by reflections.

\begin{lem}\label{lem:tower}
Given a configuration $T$ in a vertical $2\times2\times\Z$ box,
if all the central horizontal squares are present as normal faces in~$T$,
then $T$ is a tower.
\end{lem}
\begin{proof}
Each layer of $T$ is one of the six local configurations
around a normal face.  But each one of these can stack vertically
only to its mirror image.  So $T$ must be generated by reflections
from any of its layers.
\end{proof}

\vvp
\section{Screw vertices in cubic polyhedra}

Because monkey-saddles do not have flanges, the flanges of
any screw vertex in a cubic polyhedron must connect
to further screws along the axis.
Therefore, any screw vertex lies in a bi-infinite \ddefine{column}
of screws with a common axis line.  Such a column $C_\sigma$
is specified by a sequence $\sigma: \Z\to\{\Scrw,\Zcrw\}$ specifying
the handedness of each vertex.  (Shifting or reversing the sequence
results in a directly congruent column; interchanging $\Scrw$ and $\Zcrw$
results in a reflected column.)

\begin{lem}
There are uncountably many cubic polyhedra with all screw vertices.
\end{lem}

\begin{figure}[ht]\centering
\begin{minipage}[b]{.27\linewidth}\centering\epsfig{figure=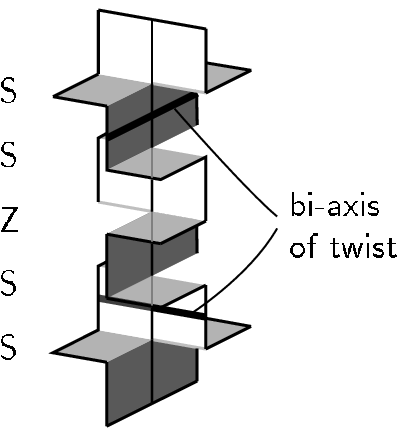}\end{minipage}
\begin{minipage}[b]{.25\linewidth}\centering{\epsfig{file=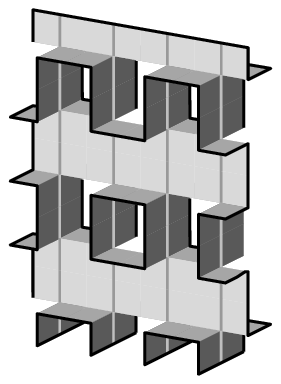}}\end{minipage}
\begin{minipage}[b]{.10\linewidth}\centering{\epsfig{file=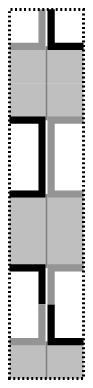}}\end{minipage}
\begin{minipage}[b]{.25\linewidth}\centering\epsfig{figure=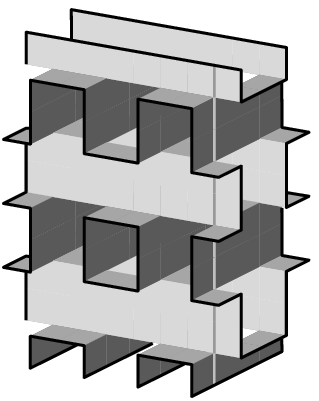}\end{minipage}
\ccaption{The column $C_\sigma$ of screw vertices corresponding the the sequence
$\sigma=\cdots\Scrw\Scrw\Zcrw\Scrw\Scrw\cdots$ (left); the slab $S_\sigma$
it generates; a planar diagram for the column (or for the slab); and the
cubic polyhedron $P_\sigma$ it generates by reflection (right).}\label{column}
\end{figure}

\begin{proof}
For any sequence $\sigma$, the column $C_\sigma$ is a configuration
in a $1\times1\times\Z$ box, which can be reflected to a complete
cubic polyhedron $P_\sigma$ by Lemma~\ref{lem:refl}.
(See figure~\ref{column}.)
These polyhedra are congruent only when the corresponding columns are.
\end{proof}

When two adjacent vertices in a column have the same handedness,
we say the column has a \ddefine{twist} along the flange joining them;
the \ddefine{normal} of a twist is the normal direction to that flange
(and the \ddefine{bi-axis} is the direction perpendicular to
both the axis and the normal).

Let $\sigma_1$ be the alternating sequence
$\cdots\Scrw\Zcrw\Scrw\Zcrw\cdots$.  Then $C_{\sigma_1}$ is the
unique untwisted column, and $P_{\sigma_1}$ is
the polyhedron $P_1$ we saw in the introduction.

Note that each $P_\sigma$ divides $\Et$ into two congruent regions. 
To see this, shift any $P_\sigma$ by $\langle 1,1,0 \rangle$;
this interchanges the components of the complement of $P_\sigma$,
but leaves $P_\sigma$ invariant.  (The polyhedron $P_0$ also
divides space into two congruent regions, as seen by a body diagonal
$\langle 1,1,1 \rangle$ translation.)

If we (repeatedly) reflect a column $C_\sigma$ in one coordinate direction
but not the other, we get a $1\times\Z\times\Z$ box which we call
the \ddefine{slab} $S_\sigma$; see figure~\ref{column}.
(We can also fill out a $1\times\Z\times\Z$ box with reflected monkey-saddles;
we call this configuration a \ddefine{sheet}.)

The \ddefine{axis} of a slab is the direction
of the axis line in any of its columns, and its \ddefine{normal} is
its direction of finite extent.  (Again, the \ddefine{bi-axis} is the
direction perpendicular to both the axis and the normal.)
Note that if two slabs in a cubic polyhedron intersect, they do
so in a common column, and thus they have the same axis direction.
Any slab is invariant under translation by two units along its bi-axis,
because it was generated by reflections in planes a unit distance apart.

\begin{figure}[ht]
\centering\epsfig{file=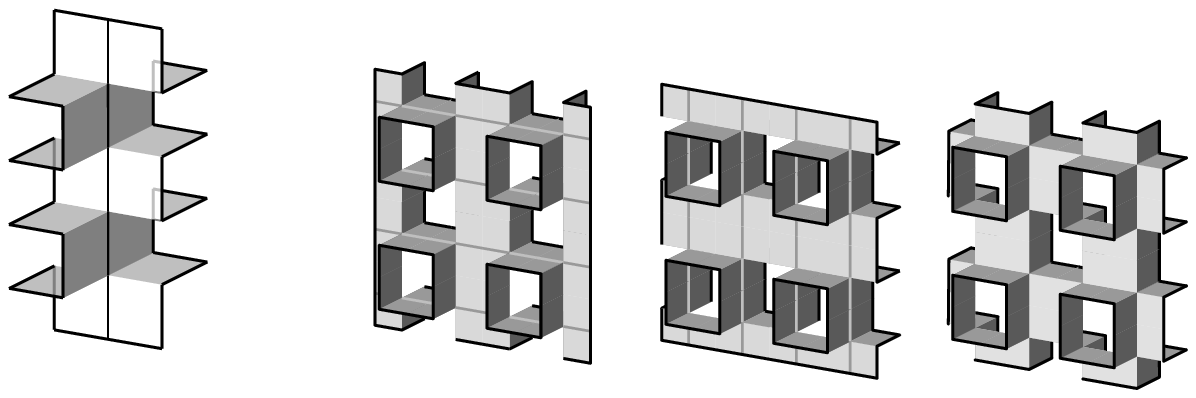}
\ccaption{The untwisted column $C_{\sigma_1}$ (left).
Any sequence of sheets of monkey-saddles (right) or
untwisted slabs $S_{\sigma_1}$ (center), which can be rotated by a quarter-turn
so that their normals agree but their axes are perpendicular,
can fit together front-to-back to form a polyhedron $P_\tau$.}
\label{ptau}
\end{figure}

We now construct a second uncountable family of cubic polyhedra,
indexed by a bi-infinite ternary sequence $\tau:\Z\rightarrow\{0,x,y\}$.
The polyhedron $P_\tau$ is built from layers in horizontal planes.
If $\tau(n)=0$, then there is a sheet of monkey-saddles in
the plane $z=n$; otherwise there is an untwisted slab $S_{\sigma_1}$
oriented to have normal $z$ and axis $\tau(n)$.
(Because $S_{\sigma_1}$ is untwisted,
it fits against a copy of itself with or without a $90^\circ$ rotation,
or against a sheet, as illustrated in figure~\ref{ptau}.)
 
\begin{lem}\label{lem:parallel}
Two adjacent parallel columns must be reflections of each other. 
\end{lem}
\begin{proof}
If two adjacent screw vertices have the same handedness,
then either they share an axis in a common column,
or they have perpendicular axes.  In two adjacent parallel
columns, corresponding vertices thus have opposite handedness, and
are mirror images by Lemma~\ref{lem:uniq}.
\end{proof}

\begin{lem}\label{lem:contig}
Three mutually perpendicular untwisted columns cannot be mutually adjacent.
\end{lem}
\begin{proof}
Two adjacent perpendicular columns touch in screws of the same handedness,
but handedness alternates along any untwisted column.
\end{proof}

\begin{lem}\label{lem:twistslab}
Let $C$ be a column of screws in a cubic polyhedron $P$.
If $C$ has a twist, then $C$ lies in a slab with the same normal as that twist. 
Moreover, if $C$ has a pair of twists separated by an odd distance,
$P$ is congruent to some $P_\sigma$, with all screw vertices.
\end{lem}

\begin{figure}[ht]
\centering{\epsfig{file=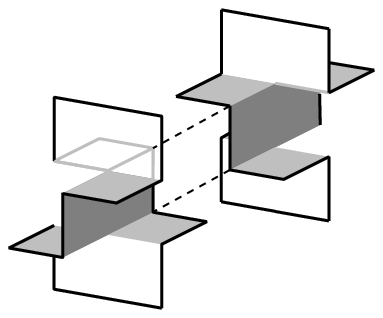}}
\ccaption{The faces on the flange of a twist in one column
are in the configuration described by Lemma~\ref{lem:crease}, and
thus must connect to another twist in a mirror image column.}\label{neartwist}
\end{figure}

\begin{proof}
Suppose there are two successive screws of the same handedness
in a column $C$.  The edge connecting them is their common axis,
a flange between two coplanar faces $f_1$ and $f_2$.
(See figure~\ref{neartwist}.)
But each $f_i$ is then of the type described by Lemma~\ref{lem:crease},
so its opposite edge is similarly a twist in a column parallel to $C$.

If a column $C_\sigma$ has twists with odd separation distance,
then their normals are perpendicular.  Applying the argument
above, adjacent to $C_\sigma$ in any direction there must be
a mirror-image column.  Repeating, we see that our polyhedron
must be $P_\sigma$.
\end{proof}

Consequently, if $S_\sigma$ is a slab in a polyhedron other than 
$P_\sigma$, the sequence $\sigma$ must be formed by pairs $\Scrw\Zcrw$ and 
$\Zcrw\Scrw$; since twists correspond to consecutive vertices of the 
same handedness, this ensures all twists lie an even distance apart.

\begin{lem}\label{lem:twistedslabs}
Suppose $P$ is a cubic polyhedron containing a twisted slab $S=S_\sigma$.
Unless $P=P_\sigma$, every slab in $P$ has the same normal as $S$.
\end{lem}
\begin{proof}
Any slab $S'$ with a different normal would intersect
$S$ along a twisted column $C=C_\sigma$; let $n$ be the normal
direction to some twist in $C$.
Both slabs $S$ and $S'$ have the same axis (along $C$) but distinct normals;
let $T$ denote the one whose normal is not equal to $n$.
The consecutive columns in the slab $T$ are reflections of one another,
so they all are twisted and have normal $n$.  By Lemma~\ref{lem:twistslab},
each lies in a slab with normal $n$, so $P=P_\sigma$.
\end{proof}

\begin{lem} Let $P$ be a cubic polyhedron with all screw vertices 
such that there is some column in $P$ with a twist. Then all columns
in $P$ are parallel, and $P$ is congruent to some $P_\sigma$.\end{lem} 

\begin{proof}
Let $C$ be a column with a twist between $\Scrw$ vertices $v_1$ and $v_2$.
By Lemma~\ref{lem:twistslab},
$C$ lies in some slab with the same normal as that twist.
Consider the vertices $v'_i$ adjacent to $v_i$ in that normal direction.
If either $v'_i$ is $\Zcrw$, then it is a reflection of $v_i$,
with axis parallel to that of $C$, so we get a parallel column.
Then this column is also in a slab,
and propogating this argument, we have the desired result.

But the $v'_i$ cannot both be $\Scrw$, for then each would
have axis perpendicular to that of $v_i$.  But then $v'_1$
and $v'_2$ would be adjacent $\Scrw$ vertices with parallel axes,
contradicting Lemma~\ref{lem:parallel}.
\end{proof}

\begin{thm}\label{thm:allb} 
Any polyhedron with all screw vertices is congruent to a 
$P_\sigma$ or a $P_\tau$.
\end{thm}

\begin{proof}
Let $P$ be a cubic polyhedron with all screw vertices.
The vertices of the cubic lattice are partitioned into columns,
which clearly can have axes in at most two directions. 
We may assume each column is untwisted, for otherwise
the last lemma would apply.
If all the axes are parallel, then $P$ is congruent to some $P_\sigma$
by Lemma~\ref{lem:parallel}. 
Otherwise, we can at least partition the vertices into planes,
within which all axes are parallel.
Within each plane, the (untwisted) columns form a slab
congruent to $S_{\sigma_0}$.
The family $P_\tau$ includes all possible ways for these 
untwisted slabs to fit together, so we are done. 
\end{proof}

\begin{figure}[ht]\centering\epsfig{file=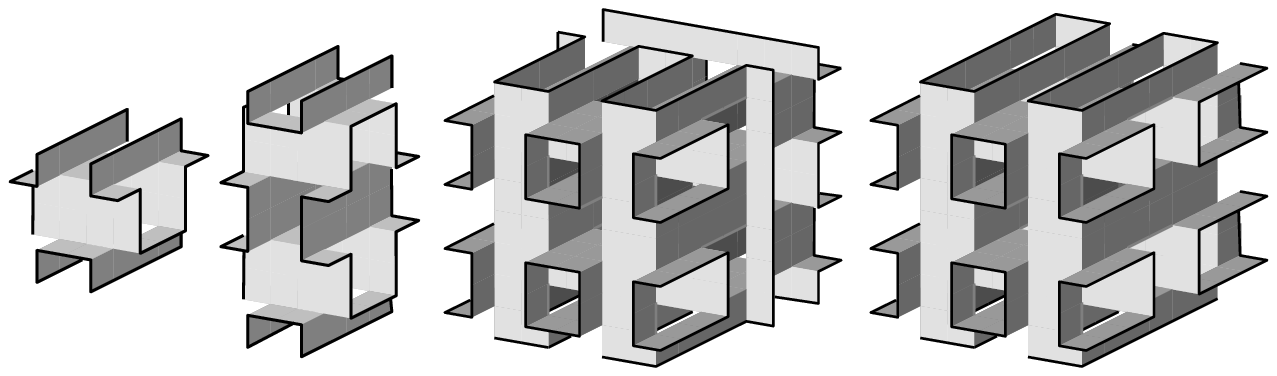}
\ccaption{The four additional new uniform polyhedra
$P_{\sigma_0}$, $P_{\sigma_2}$, $P_{\tau_1}$ and $P_{\tau_2}$.}\label{Parchs}
\end{figure}

\vp
Note that any cubic polyhedron which is uniform
in the sense of having regular faces and a symmetry acting 
transitively on the vertices
must either be $P_0$ or have all screw vertices and be
some $P_\sigma$ or $P_\tau$, where the
sequences $\sigma$ and $\tau$ must have transitive symmetry group.
The possibilities (up to isometry) are exactly
$$\sigma_0=\cdots\Scrw\Scrw\Scrw\cdots,\quad
\sigma_1=\cdots\Scrw\Zcrw\Scrw\Zcrw\cdots,\quad
\sigma_2=\cdots\Scrw\Scrw\Zcrw\Zcrw\Scrw\Scrw\cdots,$$
$$\tau_0=\cdots x x x\cdots,\quad
\tau_1=\cdots x y x y\cdots,\quad
\tau_2=\cdots x x y y x x\cdots.$$
Of course, $P_{\tau_0} = P_{\sigma_1} = P_1$ is the uniform
polyhedron already mentioned.  The four others,
$P_{\sigma_0}$, $P_{\sigma_2}$, $P_{\tau_1}$ and $P_{\tau_2}$,
are additional new uniform polyhedra, shown in figure~\ref{Parchs}.
Of course, any all-screw polyhedron in the uncountable
families $P_\sigma$ and $P_\tau$ is semi-regular in the weaker
sense of having regular faces and congruent vertex figures.

\vvp
\section {Operations on cubic polyhedra}

In order to understand cubic polyhedra that mix monkey-saddle and
screw vertices, we next define two operations.  Pushing a tower
changes monkey-saddle vertices to screw vertices (and vice versa)
within the four columns of the tower.  Inserting a slab cuts
a polyhedron along an appropriate plane, moves the two pieces
apart, and adds a new slab of screws between them.
Our \MainTheorem/ will say that
applying these operations in turn, starting from $P_0$, suffices
to create any generic cubic polyhedron. 

\vp
\subsection{Pushing towers}

\begin{figure}[ht]
\centering{\epsfig{file=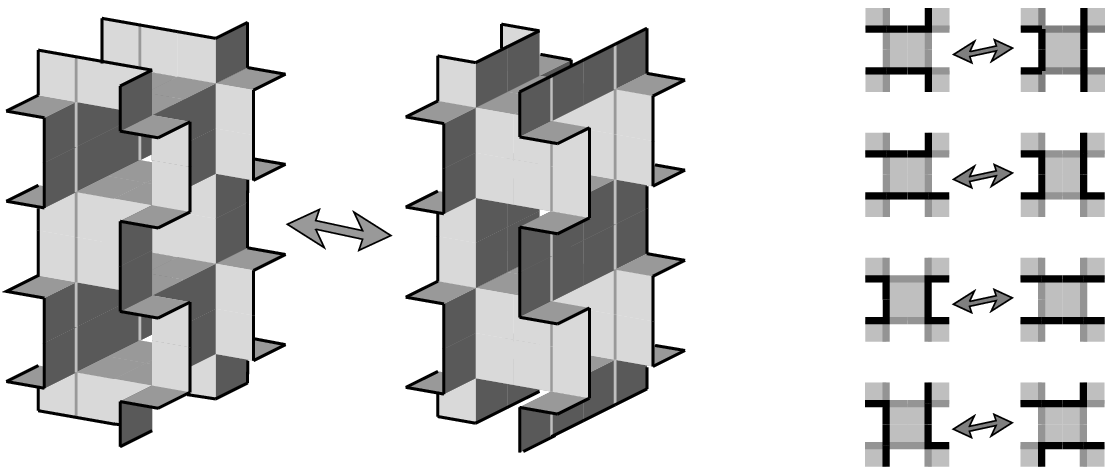}}
\ccaption{At left, we see two typical towers
(corresponding to the top pair of diagrams at right.
They differ from each other (as do the other pairs)
by a push operation, which moves all the vertical faces
one unit up (or, equivalently, down).}\label{pushtower}
\end{figure}

Remember that a tower is one of six possible configurations
in a $2\times 2\times\Z$ box obtained by mirror reflection
from the local configuration around a normal face.
Typical towers are shown in figure~\ref{pushtower}.
Of course, a screw vertex in a tower is part of a column
contained in the tower.  Note also that the vertical
faces within a tower are present or absent in checkerboard fashion.

Given a tower in a cubic polyhedron $P$, we \ddefine{push} 
the tower by moving each face within the tower one unit
along the axis of the tower.  Since all horizontal
faces are present in the configuration, we need not worry about them.
So we can equivalently describe the push
as removing all the vertical faces within the tower 
and filling in vertical faces where before there were none.
In figure~\ref{pushtower} we show the result of pushing on each
diagram from figure~\ref{towerdiag}. The following lemma should be clear:

\begin{lem} 
The result of pushing any tower is again a cubic polyhedron. 
After pushing, monkey-saddle vertices in the tower become
screw vertices and vice versa.
\qed\end{lem}

As an exercise, we note that each polyhedron in figure~\ref{lotsexamples}
can be created by applying tower pushes to $P_{0}$,
except the one in the upper right.
(That one has no towers,
and can be described instead by the techniques below.)

The columns within a tower are untwisted.  Conversely,
the following lemma shows that it is not hard to find
towers around untwisted columns.

\begin{lem}\label{col2tow}
Suppose $P$ is a cubic polyhedron, and $f$ is a normal face of $P$
at least one of whose vertices is a screw $v$ lying in an untwisted column.
Then $f$ is in a tower with axis along the normal to $f$.
\end{lem}
\begin{proof}
The local configuration around $f$ looks like one of the six configurations
in figure~\ref{towerdiag}.  The first of the six is ruled out, as it
has only $\Msad$ vertices.  Let $v$ be a screw vertex of $f$,
call its axis direction (normal to $f$) vertical, and assume without
loss of generality that $v$ is an $\Scrw$ vertex.  Above and
below $v$ are $\Zcrw$ vertices $v_\pm$, since this column has no twists.
Thus there are parallel faces $f_\pm$ in $P$ just above and below $f$.
Since the edges of $f_\pm$ incident to $v_\pm$ are both creases,
the configuration near $f_\pm$ looks like one of the diagrams in
figure~\ref{towerdiag} or~\ref{flangediag}.  But none of
the possibilities with flanges (figure~\ref{flangediag})
can fit above or below the normal face $f$.  Thus $f_\pm$ must
also be normal, and in fact mirror images of $f$.  Continuing
in this way, we find the configuration to be part of a tower.
\end{proof}

\begin{lem}\label{lem:col2tow2}
Let $P$ be any cubic polyhedron, and $C$ be any column of
screw vertices in $P$.  If $C$ does not lie in a slab or in a tower,
then $C$ lies sandwiched between two slabs whose parallel axes are
perpendicular to the axis of $C$.
\end{lem}
\begin{proof}
If the column $C$ has a twist, it lies in a slab.
So suppose $C$ is untwisted, and call its axis direction vertical.
If, among the horizontal faces incident to $C$, one is normal,
the previous lemma applies, and $C$ is part of a tower.
Thus, we may assume that every horizontal slice through $C$ is in
the configuration \inlinexb{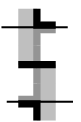}.
But then the screw vertices to either side of $C$
must lie in slabs, for they lie in columns with axes
(the thin black lines) perpendicular to that of $C$,
and such columns must occur at every level of $C$.
Therefore the column $C$ lies between two slabs, as desired.
(It is possible that $C$ is also part of a slab itself.)
\end{proof}

Our \MainTheorem/ will say that any cubic polyhedron either is a $P_\sigma$ or
is obtained from a $P_\tau$ by pushing some towers and then
inserting slabs.

Let us consider how to apply tower pushes to a $P_\tau$.
The polyhedron $P_0$ has a checkerboard of 
towers in each coordinate direction.
Any other $P_\tau$ has slabs with horizontal axis,
and thus has no vertical towers.
We find horizontal towers at any level
except where $\tau$ alternates from $x$ to $y$.
The set of horizontal towers is always nonempty, except in $P_{\tau_1}$.

If we are starting with $P_0$, we first
push any (infinite, finite or empty) subset of the vertical towers.
(A nontrivial push will create some flanges and thus destroy some
horizontal towers, but it may also create new ones.)

We now take the resulting polyhedron (or our starting $P_\tau$)
and push any subset of its towers in the $y$-direction.
(Note that if we push a tower at a level where $\tau$ had a $0y$ or $yy$,
we will destroy the slab at that level.)
Finally, we push any subset of the remaining towers in the $x$-direction.

We call any cubic polyhedron resulting from this
procedure a \ddefine{\pushed/}.

\vp
\subsection{Removing and inserting slabs}\label{sectslab}

\begin{figure}[ht]
\centering{\epsfig{file=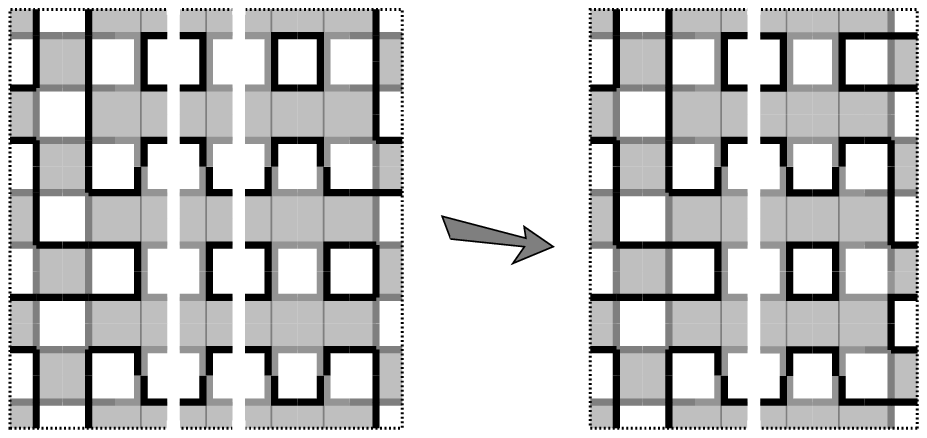}}
\ccaption{To remove a slab, we delete it and join the
two resulting half-spaces together by a translation in the
direction of the slab normal plus the bi-axis.  In this diagram
(in the axis/normal plane) the bi-axis translation appears as
swapping all edge colors in the right half-space.}\label{removeslab}
\end{figure}

Given a cubic polyhedron $P$, we \ddefine{remove} a slab by 
the following operation:
delete the slab completely, and then join the two (half-space)
components of the complement together by translating one of them
by one unit along the slab normal (to bring them together)
plus one unit along the slab bi-axis (to let them match).
For example, in figure 
\ref{removeslab} we show the result of removing a slab from the
cubic polyhedron partly indicated by the planar diagram.
It is not hard, using such diagrams, to prove the following:

\begin{lem}\label{remslab}
Let $P$ be a cubic polyhedron containing a slab $S$;
the result of removing $S$ is again a (well-defined) cubic polyhedron.
\qed\end{lem}

We define \ddefine{inserting} a slab to be the inverse of this operation.
After the next lemma, we will explain exactly what we mean by
inserting an infinite sequence of slabs, but we note now that
any $P_\tau$ can be obtained from $P_0$ by inserting
(perhaps infintely many) untwisted slabs.

We will now examine the conditions under which a slab can be inserted
into a \pushed/, and will find that the slab to be inserted is usually
uniquely determined.

\begin{lem}\label{lem:slabpusher}
Let $P$ be a \pushed/.
If $P$ has columns in all three directions, no slab can be inserted in $P$.
If $P$ has no vertical columns, then slabs can be 
inserted at any horizontal plane $z=i+\frac{1}{2}$, $i\in\Z$.
\end{lem} 
\begin{figure}[ht]
\centering{\epsfig{file=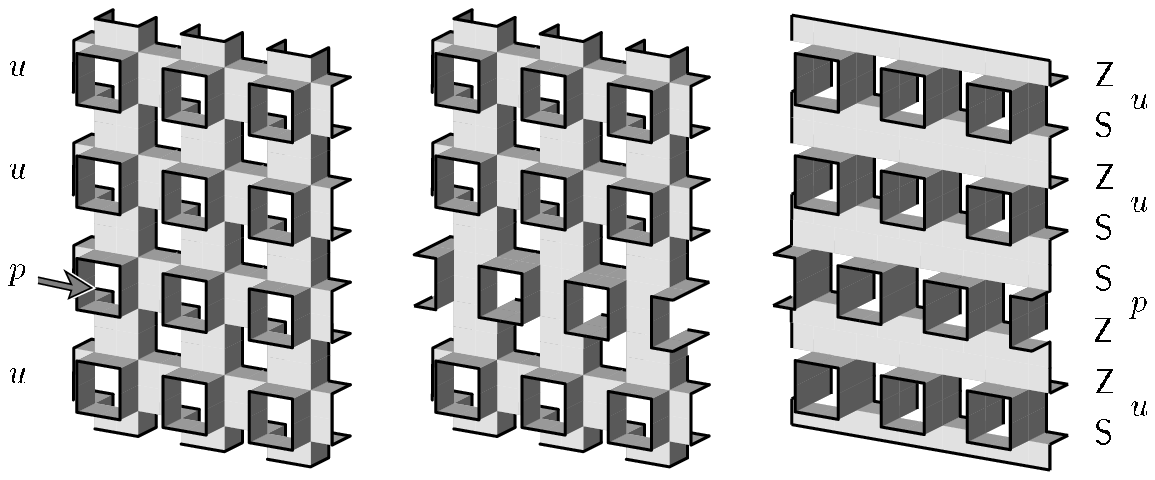}}
\ccaption{The sequence of tower pushes in a plane determine which 
slabs can be inserted. We start with the left figure, exhibiting the trivial
pattern of squares on the front boundary plane of the image.
We push towers with axes in this plane, producing the middle figure.
We will then be able to attach the slab illustrated at right,
since they both show the same pattern of squares.}
\label{slabpush}
\end{figure}

\begin{proof}
If $P$ has columns in all three directions, then any plane
is cut perpendicularly by the axis of some column, and thus cannot be
a candidate for slab insertion.

If there are no vertical columns, consider a horizontal plane
$\plh$ given by $z=i+\frac{1}{2}$.
We will refer to the intersection of $P$ with $\plh$ as the \ddefine{pattern}
of $P$ in $\plh$.  This shows us the set of faces in $P$ which are
bisected by $\plh$, and is always a union of squares.
If we insert a slab along $\plh$,
the boundary of the slab must have the same pattern.

If $P$ is a $P_\tau$, the pattern we see is always the \ddefine{trivial pattern}
consisting of a square array of squares, as in figure~\ref{slabpush} (left).
(This is why the three possible layers in a $P_\tau$ can fit together
in any order.)

To determine the pattern of $P$ in $\plh$, we need only look
at which tower pushes in the plane $\plh$ have been performed.
(Note that horizontal tower pushes which occur one layer higher or lower
in $P$ will change the vertices just above and below $\plh$, but won't
affect the faces cut by $\plh$ or thus what slab we can insert.
Also, depending on $\tau$, pushes in $\plh$ may or may not be possible.)

If any tower within $\plh$ has been pushed, then all tower pushes
at that level (or at levels just above and below $\plh$) must
have parallel axes.  The only exception is that if all $x$-towers
within $\plh$ are pushed, then we again see the trivial pattern in $\plh$
(shifted by one unit), and any subset of the (new) $y$-towers in $\plh$ can
be pushed.

If the pattern of $P$ in $\plh$ is trivial (meaning that either no
towers or all towers in $\plh$ have been pushed),
we can insert an untwisted slab with axis $x$ or $y$.
Indeed, we can insert any finite number of 
untwisted slabs, indexed by a finite sequence of $x$'s and 
$y$'s, or we can delete the half-space of $P$ to one side of $\plh$
and insert a half-space of some $P_\tau$ corresponding to
an infinite sequence of $x$'s and $y$'s.  Of course, the result
of such an insertion is merely another \pushed/.

On the other hand, if the pattern is nontrivial, then
all tower pushes have occurred with parallel axes, say along $x$,
and are at an even distance apart from one another within $\plh$.
We now define a bi-infinite sequence $\omega$ of letters $u,p$:
reading along the $y$-direction in $\plh$,
each tower with axis in the $x$-direction was either $u$npushed or $p$ushed.
The substitution rules $u\mapsto\Zcrw\Scrw$, $p\mapsto\Scrw\Zcrw$
convert $\omega$ to a sequence $\sigma$.  Then the unique slab
that can be inserted along $\plh$ is $S_\sigma$,
as illustrated in figure~\ref{slabpush}.
Moreover, we can insert any finite number of $S_\sigma$'s, each a
reflection of the next, or we can delete the half-space to one side of
$\plh$ and insert a half-space of $P_\sigma$.
\end{proof}

If $P$ is a \pushed/ with no vertical columns, we can now explain
what we mean by insertion of a (possibly infinite) set of slabs.
Beginning near the origin and working outwards, we examine the planes 
$z=\pm (i-\frac{1}{2})$, for $i=1,2,\ldots$.
If a given plane has the trivial pattern,
we can insert any finite or infinite sequence of 
untwisted slabs, indexed by choice of axis, $x$ or $y$.
Otherwise, there is a sequence $\sigma$ associated with the plane,
as described in the proof of Lemma~\ref{lem:slabpusher},
and we can insert any finite or infinite
number of copies of $S_\sigma$.
(To insert an infinite sequence of slabs at one level,
we throw out the half-space bounded by the plane away from the origin,
and then stop the procedure on that side.)
In our \MainTheorem/, we will show that any cubic polyhedron
can be obtained from a \pushed/ in this way.

\vvp
\section{Classifying Cubic Polyhedra} 

We now show that, given any cubic polyhedron, we can
apply the operations of the previous section in reverse
to simplify it.
First we consider cases where tower pushes alone can bring us back to $P_0$.

\begin{lem}\label{lem:push1}
Let $P$ be a cubic polyhedron in which all columns (if there are any)
are untwisted with vertical axis.  Then $P$ can be obtained
from $P_0$ by pushing some set of vertical towers.
\end{lem}
\begin{proof}
The polyhedron $P$ must be generated by reflection from
any $\Z\times\Z\times1$ box around a horizontal plane $\pli$,
since the screws lie in untwisted vertical columns and
the monkey-saddles also lie in infinite vertical lines.  It follows that
every horizontal face in $P$ is normal and part of a vertical tower.

We will assign a \ddefine{parity} modulo $2$ to each tower as follows.
Arbitrarily assign parity~$0$ to one tower, centered at some point $(x,y)$.
We may assume that the planar diagram in $\pli$ at this tower
is~\inliness{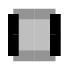}. Any other tower is centered at
some $(x+a, y+b)$ with $a\equiv b \pmod 2$.
To this tower we assign parity $a$
if it has the same diagram~\inliness{towdiag1}, and parity $a+1$
if it has the other diagram~\inliness{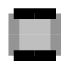}.
Note that two adjacent towers intersect in a column of screws
if they have opposite parity, as in~\inliness{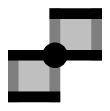},
while they intersect in a line of monkey-saddles if they have
the same parity, as in~\inlines{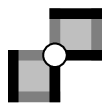}.

Pushing a tower switches its parity and leaves all other towers unaffected.
So we now simply push all towers of parity~$1$;
this leaves all towers with the same parity,
so we now have all monkey-saddle vertices, as desired.
\end{proof}

\begin{lem}\label{lem:push3}
Let $P$ be a cubic polyhedron with columns in all three coordinate directions.
Then we can push some set of vertical towers to eliminate all
vertical columns.  All columns in the resulting polyhedron $P'$ are untwisted
and horizontal, and $P'$ has at least one vertical tower (of monkey-saddles).
\end{lem}

\begin{proof}
There can be no slabs in $P$, and thus all columns are untwisted.
We will first find pushes
to remove all vertical columns.

Choose any horizontal plane $\pli$ in the cubic lattice,
and consider the set $V$ of all screw vertices in $P$ with horizontal axis.
Since these lie in horizontal columns, the projection of $V$ to $\pli$
consists of infinite lines of vertices, including at least one line in
each direction.
Its complement (within the set of all vertices in $\pli$)
is a (possibly infinite) union of rectangular regions.
Each region is finite or
singly infinite in each direction, and thus has connected boundary.

In figure~\ref{dotslice} we illustrate a schematic of the slice in $\pli$.
Somewhere above or below each occurence of the symbol~\inliness{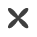}
there is a screw vertex with horizontal axis.
At the locations marked~\inline{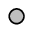} there is either a monkey-saddle 
vertex~\inline{adot} or a screw vertex with vertical axis~\inline{bdot}. 
\begin{figure}[ht]
\centerline{\epsfig{file=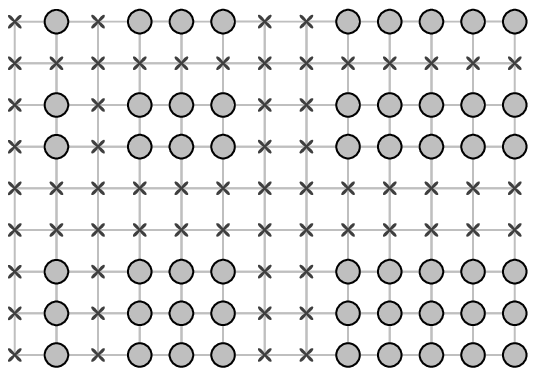}}
\ccaption{A schematic for a typical slice through a cubic polyhedron with 
columns in three directions.}\label{dotslice}\end{figure}

Now consider some $\alpha\times\beta$ rectangle $R$ of vertices
in the complement of the projection of $V$, and assume $\alpha\ge\beta$.
The configuration around $R$ is a $\alpha\times\beta\times1$ box
which necessarily appears reflected in $P$ to fill out a
$\alpha\times\beta\times\Z$ box.  (Indeed, any screw
vertex in $R$ is part of an untwisted vertical column, and any monkey-saddle
vertex lies above and below further monkey-saddles by the definition of $R$.)

We will prove that we can push vertical towers within $R$ to
eliminate all screw vertices in $R$.
Within $R$, horizontal faces appear in a checkerboard pattern, as
seen in figure~\ref{pushregion}.  Those along the boundaries
of~$R$ (shown in lighter gray in the figure) are not in towers,
since they are incident, at some level, to horizontal screws in~$V$.
But all interior horizontal faces (in darker gray)
lie in vertical towers, by Lemma~\ref{col2tow} and the definition of~$R$.

Since $P$ has no slabs, by Lemma~\ref{lem:col2tow2}
any column of screws must be in a tower.
If $\beta=1$ then there are no towers in $R$, so there can be
no screws in $R$, and we are done.

We now assume $\beta>2$; we will return later to the case $\beta=2$.
Note that interior vertices of $R$ are incident to exactly two towers, while
vertices along an edge are incident to exactly one.  A corner vertex
in $R$ is incident to either one tower or zero towers, depending
on whether or not there is a horizontal face in that corner of $R$.

We assign a parity modulo $2$ to each tower in $R$ exactly
as in the proof of Lemma~\ref{lem:push1},
and then push all towers of parity~$1$. 
Again, this leaves all towers within $R$ with the same parity,
meaning that the interior of $R$ has only monkey-saddle vertices.
(Beginning with any polyhedron represented by the first configuration
in figure~\ref{pushregion}, we have now obtained the second one.)
\begin{figure}[ht]\centering
\begin{minipage}[b]{.2\linewidth} 
  \centering{\epsfig{file=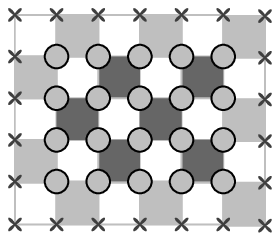}}\end{minipage}\hfill
\begin{minipage}[b]{.2\linewidth} 
  \centering{\epsfig{file=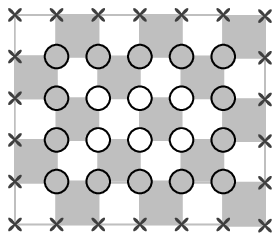}}\end{minipage}\hfill
\begin{minipage}[b]{.2\linewidth} 
  \centering{\epsfig{file=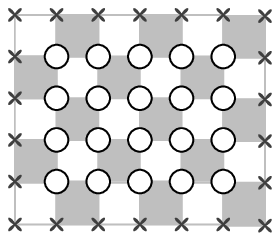}}\end{minipage}\hfill
\begin{minipage}[b]{.2\linewidth} 
  \centering{\epsfig{file=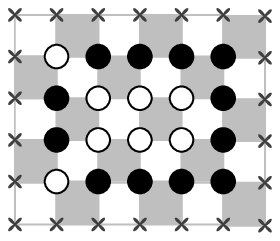}}\end{minipage}
\ccaption{Four configurations illustrating the effect of
pushing towers.}\label{pushregion}
\end{figure}

If all boundary vertices of $R$ are also now monkey-saddles
(as in the third configuration in figure~\ref{pushregion}),
we have proved the claim, having eliminated all vertical columns
from this portion of the cubic polyhedron.

Otherwise, we claim that each screw vertex in $R$ is incident to exactly
one tower, while each monkey-saddle is incident to zero or two towers.
(This is the case in the last configuration in figure~\ref{pushregion}.)
Then pushing all towers in $R$ will remove all vertical columns, as desired.
(In other words, if we originally had pushed the towers of parity~$0$
instead of those of parity~$1$, we would have obtained all monkey-saddles.)
The claim is equivalent to saying that there are screws along the
edges of $R$, and screws in corners where a horizontal face is present,
but monkey-saddles in the other corners.

Suppose there is a screw in an inappropriate corner, that is,
we see the configuration~\inliness{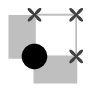}.  
At some levels above or below $\pli$, there are horizontal
columns along one edge of $R$; at other levels there are perpendicular
horizontal columns along the other edge of $R$.
But these two possibilities cannot happen at successive levels,
since then we would have
three mutually perpendicular and adjacent untwisted columns,
contradicting Lemma~\ref{lem:contig}.
This means that at some intermediate levels, we have
monkey-saddles at all three positions marked~\inliness{abovebdot}.
This means we see the configuration $F:=$\inliness{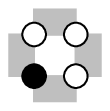},
which is impossible as we noted in Section~\ref{config}.

Remember we have assumed there is at least one screw on the boundary of $R$.
We cannot have the configuration~\inliness{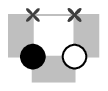}, since, again,
at some level we will have monkey-saddle vertices at the positions
marked~\inliness{abovebdot} and thus would have the forbidden 
configuration $F$. Similarly, we cannot have 
the configuration~\inliness{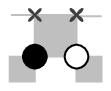}, for there are only monkey-saddles
in the (non-empty) interior of $R$,
and again we find the forbidden configuration $F$.
Thus, along an edge of $R$ with one screw vertex,
we continue to see screws until we reach a corner.

Finally, we check that we cannot have a monkey-saddle at an
inappropriate corner: we cannot see the configuration~\inlinex{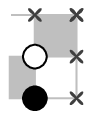},
for again at some level the positions marked~\inlines{abovebdot}
must be monkey-saddles, giving $F$ yet again.

So all corners must have the appropriate vertices,
either~\inliness{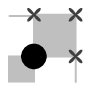} or~\inliness{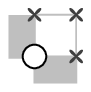}.
In the latter case, we cannot have the configuration~\inlines{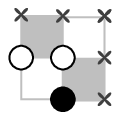},
since we have only monkey-saddles
in the interior of $R$, giving $F$ once again.
Therefore, since the boundary of $R$ is connected and
contains at least one screw,
it consists entirely of screws
(except monkey-saddles in the appropriate corners).

Finally, suppose $\beta=2$.  At the end of $R$, we must see one of 
the configurations at the left in figure~\ref{fig:case2xa}. 
Note that the configuration $F':=$~\inliness{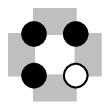} is forbidden
just like $F$.
Since we cannot have either $F$ or $F'$, by induction the towers in the 
interior of $R$ must consist of a sequence of the two configurations 
in the middle in figure~\ref{fig:case2xa}, and the vertices of each tower 
are either all screws or all monkey-saddles. Moreover, each 
vertex in $R$ belongs to just one tower. We simply push the towers 
that consist of screw vertices and leave the others alone. We are 
left with all monkey-saddles in $R$.
\begin{figure}[ht]
\centering\epsfig{file=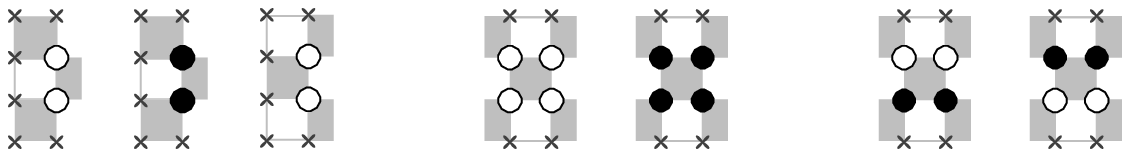}
\ccaption{At left and middle, configurations 
arising when $R$ is $\alpha\times2$; at middle and right, 
configurations arising when $R$ is a bi-infinite strip of width two,
as in the proof of Lemma~\ref{lem:push2}.}
\label{fig:case2xa}
\end{figure}

If we repeat the procedure described so far for each rectangle $R$,
we can push towers to eliminate all screws with vertical axis, and
obtain the polyhedron $P'$.  It contains no vertical columns, but does
contain the vertical towers we pushed in order to obtain it from $P$.
\end{proof}

\begin{lem}\label{lem:push2}
Let $P$ is a cubic polyhedron with no vertical columns and no slabs.
Then $P$ is a \pushed/.  If $P$ includes a vertical tower (of monkey-saddles),
then it is actually a \pushedz/.
\end{lem}

\begin{proof}
If $P$ has columns only in a single direction,
we can apply Lemma~\ref{lem:push1}.
So suppose $P$ has columns along both $x$ and $y$.
The proof will be very similar to the previous one.

Let $\pli$ be the $xz$-plane, and $V$ the set of screws
with axis in the $x$-direction.  Projecting $V$ to $\pli$,
we now find that the complement is a union of horizontal
$\Z\times\beta$ strips $R$.  When $\beta=1$, since $P$ has
no slabs, Lemma~\ref{lem:col2tow2} again guarantees that
the region $R$ contains no screws.  So we may assume $\beta>1$.

First, we will suppose that $P$ does include a vertical tower.
This tower gives a vertical strip $S$ (of width $2$)
in which we see only monkey-saddles.  The complement of $S$ in each $R$ is
two semi-infinite strips $R_\pm$.  Each of these regions has corners
and connected boundary, as in the previous proof.

The argument then proceeds as before, with the monkey saddles
along $S$ playing the same role as the monkey-saddles we earlier
found at some level along each boundary of $R$.
We assign parities to the towers in $R_\pm$, and push the towers of parity~$1$.
Then if any screws remain along the boundary, we push all towers.
The same arguments guarantee that this removes all columns from $R_\pm$.

Repeating for each region $R_\pm$, we eventually remove all columns
along $y$.  We can then use Lemma~\ref{lem:push1} to remove the
remaining columns (along $x$).

If $P$ contains no vertical tower, the argument is a bit more difficult.
We work directly on the strips $R$ of height $\beta\ge2$.
First assume $\beta>2$.  Assign parities to the towers in $R$,
and push those of parity~$1$.  Fix one boundary component of $R$;
the usual arguments imply that along this boundary now
either all vertices are screws or all are monkey-saddles.
If both boundary components have screws, push all towers in $R$
to eliminate all the screws.  However, if one boundary has
screws and the other has monkey-saddles, there is no way to
push towers to get rid of all screws; we are left with a slab
(untwisted, with axis $y$ and normal $z$) along this boundary of $R$.

If $\beta=2$, the situation is similar.
Again, since we cannot have the configurations $F$ or $F'$, 
the towers in $R$ either are a sequence of the two 
configurations in figure~\ref{fig:case2xa} (center) or are a 
sequence of the two configurations figure~\ref{fig:case2xa} (right).
In the first case, we push the towers 
consisting of all screw vertices as before, removing all columns from $R$.
In the second case, we push all copies of the far-right configuration.
(In either case, the push we use is in fact the usual push of all
towers of parity~$1$, although the accounting is more difficult here.)

In summary, within in each region $R$, we push towers to eliminate
all columns with axis $y$ except, perhaps, for a single untwisted
slab along one boundary of $R$.  We arrange (by pushing all towers
in $R$ if necessary) that these slabs are always at the top
of each region $R$.  It follows that these slabs are separated by
distance at least $3$ from each other (since each $R$ was at
least $2$ units high, and they are separated by at least one layer).
Call the resulting cubic polyhedron $P'$.

Because of these slabs, we cannot now apply Lemma~\ref{lem:push1}
to $P'$.  We could remove the slabs and then apply that lemma.
This would show our original $P$ could be obtained from $P_0$
by a combination of tower pushes and slab insertions, but would
not demonstrate that it is a \pushed/.

Instead, we repeat the argument we have just given, but with the
roles of $x$ and $y$ interchanged, starting from $P'$.
We must be careful because $P'$ does contain slabs.  But we get
$\beta\ge2$ for each strip $R$ directly, since the slabs (the
only screws with axis $y$) are separated vertically by $3$ units.

We get rid of internal columns in $R$ by pushing towers of parity~$1$
as usual.  Along a boundary of $R$, we again claim that either all
vertices are screws or all are monkey-saddles.  We can no longer
find some height above $\pli$ at which there are monkey-saddles
just outside $R$, since $R$ is now bounded by slabs.
\jmp{Here's where we might wish we had listed
all the legal configurations around missing square.}
But these slabs are untwisted so the configuration~\inliness{forbid7}
is still forbidden, and the argument goes through as before.

We can thus find tower pushes in the $x$ direction taking $P'$
to a polyhedron with no columns except untwisted horizontal slabs.
This is a $P_\tau$.  So $P'$, and then our original $P$, is a \pushed/.
\end{proof}

\begin{cor}\label{cor:push}
Let $P$ be a cubic polyhedron with no slabs.  Then $P$ is a \pushed/.
\end{cor}
\begin{proof}
The columns of $P$ must be untwisted.  If they have parallel axes,
apply Lemma~\ref{lem:push1}.  If they occur in two directions,
apply Lemma~\ref{lem:push2}.  If they occur in all three directions,
apply Lemma~\ref{lem:push3} and then Lemma~\ref{lem:push2}.
\end{proof}

\begin{mainthm}\label{mainthm}
Any cubic polyhedron $P$ is a $P_\sigma$ or is obtained by
inserting slabs into some \pushed/.
\end{mainthm}

\begin{proof}
If a cubic polyhedron $P$ has all screw vertices, then it is
a $P_\sigma$ or a $P_\tau$, and we are done.

Otherwise, the slabs in $P$ do not fill out all of $P$.
Suppose there are two slabs in $P$ with different normals.
They intersect along a column which we will call vertical.
All columns must be vertical, since a horizontal column would
cut one of the slabs.  By Lemma~\ref{lem:twistedslabs}
all columns are also untwisted (or we would have a $P_\sigma$).
Thus Lemma~\ref{lem:push1} applies to show $P$ is a \pushed/.

Thus we may assume all slabs in $P$ have the same normal.
Call this normal direction vertical and the slabs horizontal.
We can remove all slabs; we are left with an $\alpha\times\Z\times\Z$
box $P''$, with $\alpha$ being bi-infinite, infinite or finite depending
on whether there were originally zero, one or two half-spaces of slabs.
We reflect $P''$ if necessary to get a complete cubic polyhedron $P'$
with no slabs, from which $P$ can be obtained by inserting slabs.
Now Corollary~\ref{cor:push} applies to $P'$, showing $P'$ is a \pushed/,
as desired. 
\end{proof}

\vvp
\section{Discrete Minimal Surfaces in the Cubic Lattice}
As we mentioned in the introduction, a complete surface built from faces
of the cubic lattice is discrete minimal if and only if each vertex of
the lattice has one of five possible configurations of faces.
We call these $\Msad$, $\Scrw$,
$\Zcrw$, $\Flat$ and $\Emty$, where $\Msad$, $\Scrw$ and $\Zcrw$ are the
vertex configurations of cubic polyhedra, $\Flat$ represents a flat
vertex with four coplanar squares, and $\Emty$ is an empty vertex,
not in the polyhedron.  It would be interesting to extend our
classification theorem to give a complete list of all such discrete
minimal surfaces; we give a few partial results here.

First, let us give some examples.  A trivial family, indexed by subsets
of $\Z$, contains any collection of horizontal planes.  These surfaces
are not connected, and do not include every vertex in the cubic lattice.

\begin{figure}[ht]\centering
 \begin{minipage}[b]{.45\linewidth}\centering\epsfig{figure=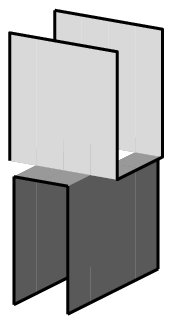,height=2in}\end{minipage}
\begin{minipage}[b]{.45\linewidth}\centering\epsfig{figure=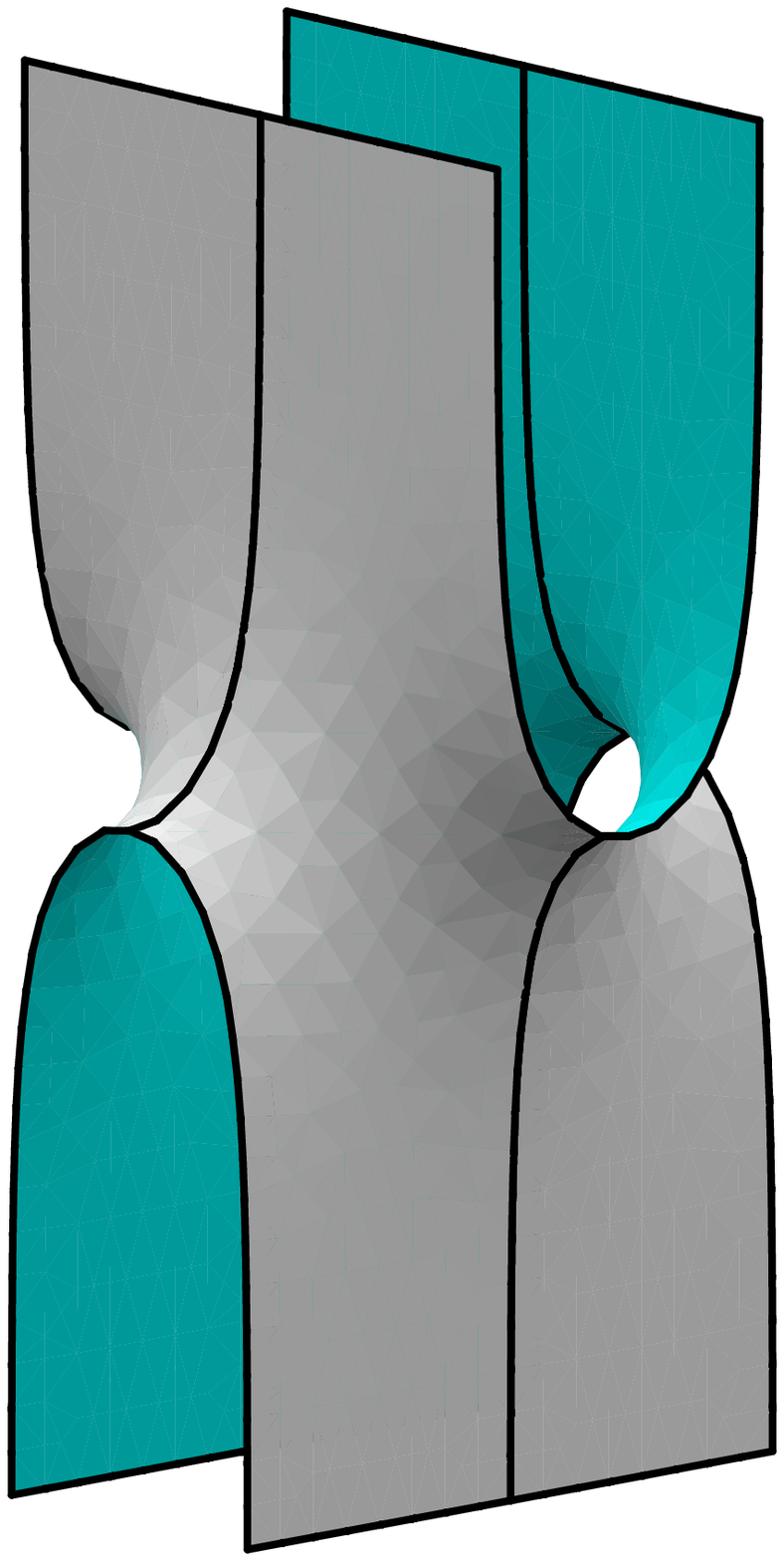,height=2in}\end{minipage}
 \ccaption{This discrete minimal surface (left) is not a cubic polyhedron, but
does relax to the singly-periodic minimal surface of Scherk (right).}
\label{scherk}
\end{figure}

We can construct a more interesting example as follows:
Attach infinite sequences of flat $\Flat$ vertices to the two flanges
of a screw vertex, and then extend this $1\times1\times\Z$ box
to a complete surface by reflection.  This is clearly a discrete
version of Scherk's doubly periodic minimal surface, and it
does relax to that surface.  (See figure~\ref{scherk}.)
In this polyhedron, there are two half-spaces containing stacked half-planes
of flat vertices, separated by a layer of screw vertices which twist
the half-planes 90 degrees.  This polyhedron is connected, and does
include every lattice vertex.

\begin{lem}
Along the normal line of a flat vertex (or any coordinate line through an empty
vertex) we find only empty vertices or flat vertices with that normal line.
\end{lem}
\begin{proof}
The empty and flat configurations are the only ones which omit
edges, and they omit edges in colinear pairs.
Thus all edges along any such line are omitted in the polyhedron.
\end{proof}

\begin{figure}[ht]\centering
 \epsfig{figure=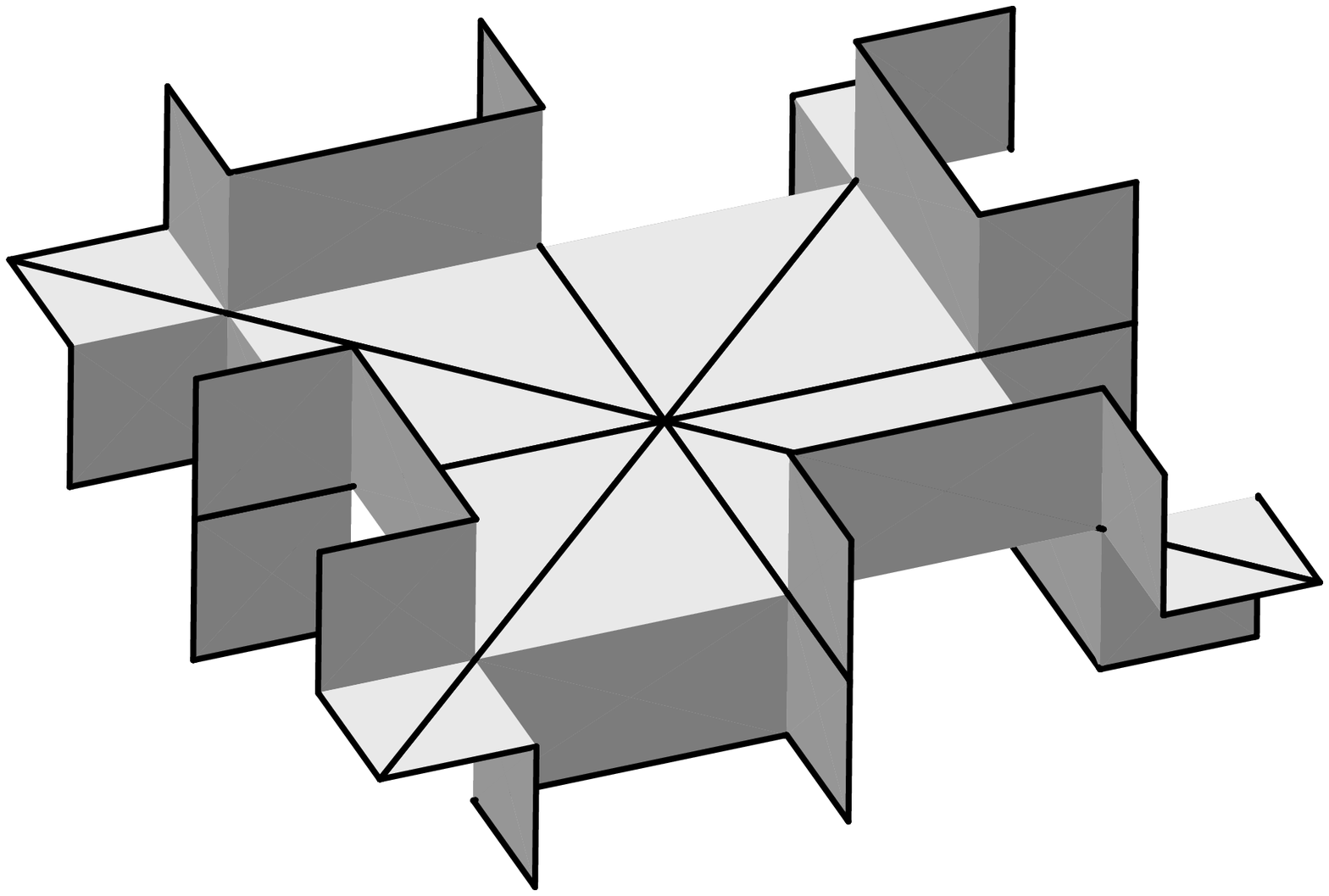,width=2in}\hspace{1in}
 \epsfig{figure=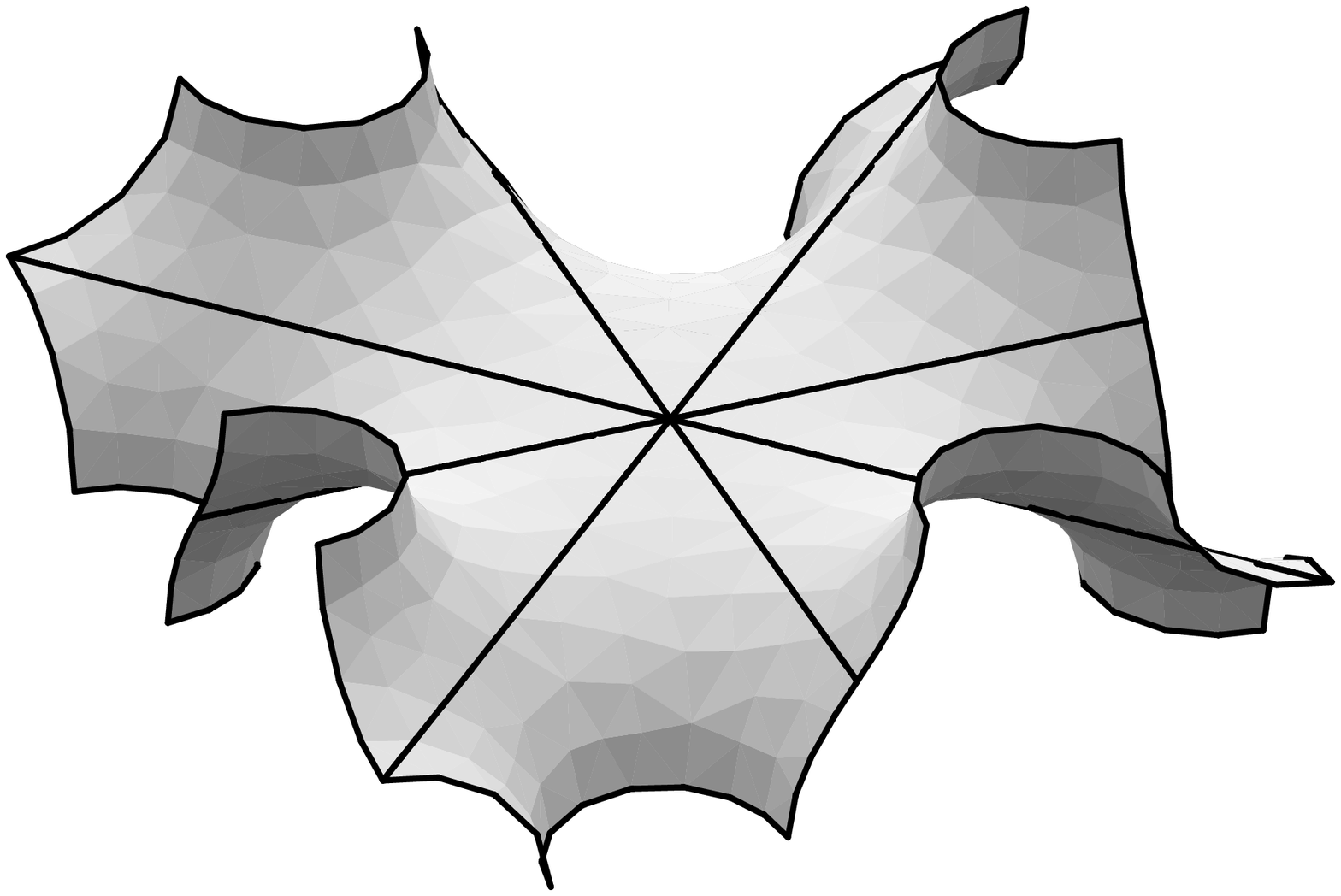,width=2in}
  \ccaption{This discrete minimal surface (left) is not a cubic polyhedron.
  Its screw vertices lie in intersecting horizontal columns.
  It has two-fold rotational symmetry around each of the eight
  lines shown through the central flat vertex, and thus four-fold rotational
  symmetry around the normal line there.  It evidently relaxes to the
  minimal surface shown (right), with the same symmetries.}
 \label{msmsxsmsm}
\end{figure}

Although this lemma constrains how flat vertices
can appear, there is still more flexibility than for
cubic polyhedra.  For instance, a flat vertex can be surrounded
on all four sides by screw vertices, as in figure~\ref{msmsxsmsm}.
This $3\times3\times1$ block can be extended by reflection.

We will not attempt here to give a complete classification of
all the discrete minimal surfaces in the cubic lattice.

\vp 

\end{document}